\documentclass{article}[12pt]
\usepackage{amsfonts,amssymb,latexsym,amsmath}
\newcommand{\mf}{\mathfrak}

\newcommand{\equi}{\Leftrightarrow}

\newcommand{\ol}{\overline}

\newcommand{\ra}{\rightarrow}
\newcommand{\Ra}{\Rightarrow}

\newcommand{\eps}{\varepsilon}

\newcommand{\mbb}{\mathbb}
\newcommand{\tn}{\textnormal}

\newtheorem{de}{Definition}[section]
\newtheorem{re}[de]{Remark}
\newtheorem{pr}[de]{Proposition} 
\newtheorem{tr}[de]{Theorem}
\newtheorem{lm}[de]{Lemma} 

\newtheorem{nt}[de]{Notation} 

\newtheorem{co}[de]{Corollary}

\newcommand{\lb}{\linebreak}
\def\vp{\rm \vspace{0.2cm}}
\def\M{\rm Max}

\def\hb{\hfill$\Box$}
\def\aut{\rm Aut} 
\def\sp{\rm Spec} 
\def\tran{\rm Trans}
\def\fl{\rm ETrans}
\def\Um{\rm Um}
\def\GL{\rm GL}
\def\SL{\rm SL}
\def\EO{\rm EO}
\def\SO{\rm SO}
\def\E{\rm E}
\def\T{\rm T}
\def\ET{\rm ET}
\def\G{\rm G}
\def\Sp{\rm Sp}

\def\k{\rm K_1}
\def\K{\rm K}
\def\T{\rm T}
\def\O{\rm O}
\def\I{\rm I}
\def\ESp{\rm ESp}
\def\EO{\rm EO}
\def\es{\rm S}
\def\sdim{\rm sdim}
\title{Local-Global Principle for Transvection Groups}
\author{A. Bak, Rabeya Basu \& Ravi A. Rao}
\date{}
\begin{document}
\maketitle
\begin{center}{\it 2000 Mathematics Subject Classification:
{13C10, 15A63, 19B10, 19B14}}
\end{center}
\begin{center}{\it Key words: 
projective, symplectic, orthogonal modules, nilpotent groups, ${\k}$}
\end{center}

\begin{abstract} In this article we extend the validity Suslin's
Local-Global Principle for the elementary transvection subgroup
of the general linear group GL$_n(R)$, the symplectic group Sp$_2n(R)$,
and the orthogonal group O$_2n(R)$, where $n > 2$, to a Local-Global
Principle for the elementary transvection subgroup of the automorphism
group Aut$(P)$ of either a projective module $P$ of global rank $> 0$ and constant local rank $> 2$, or of a nonsingular symplectic or orthogonal module $P$ of global hyperbolic rank $> 0$ and constant local hyperbolic rank $> 2$. In Suslin's results, the local and global ranks are the same, because he is concerned only with free modules. Our assumption that the global (hyperbolic) rank $> 0$ is used to define the elementary transvection subgroups. We show further that the elementary transvection subgroup ET$(P)$
is normal in Aut$(P)$, that ET$(P) =$ T$(P)$ where the latter denotes the full
transvection subgroup of Aut$(P)$, and that the unstable K$_1$-group K$_1($Aut$(P)) =$ Aut$(P)/$ET$(P) =$  Aut$(P)/$T$(P)$ is nilpotent by abelian, provided $R$ has finite stable dimension. The last result extends previous ones of Bak and Hazrat for GL$_n(R)$, Sp$_2n(R)$, and O$_2n(R)$.

An important application to the results in the current paper can be found in the 
work \cite{br}. Here the last two named authors have studied the 
decrease in the injective stabilization of classical modules over a 
non-singular affine algebra over perfect C$_1$-fields. We refer the reader 
to that article for more details. 
\end{abstract}

\section{Introduction} ~~~~In 1956, J-P. Serre asked if a finitely generated 
projective module over a polynomial ring over a field is free. 
This is known as Serre's problem on projective modules.
It was affirmatively proved by D. Quillen and A. Suslin independently in 1976. 
Now it is known as the Quillen-Suslin Theorem. Quillen established the 
following Local-Global Principle in his proof of Serre's problem in 
\cite{qu}. \vp 

{\bf Quillen's Local-Global Principle:} {\it 
A finitely presented module over a polynomial ring $R[X]$ over a commutative 
ring $R$ is extended if and only if it is locally extended over the 
localization of $R[X]$ at every maximal ideal of $R$. } \vp

We shall be concerned with the matrix theoretic version of this theorem. It 
was established by Suslin in his second proof of the Serre's problem in 
\cite{SUS}. 
\vp

{\bf Suslin's Local-Global Principle:} 
{\it Let $R$ be a commutative ring with identity and let 
$\alpha(X)\in {\GL}_n(R[X])$ with $\alpha(0)={\I}_n$. If 
$\alpha_{\mf{m}}(X)\in {\E}_n(R_{\mf{m}}[X])$, for every maximal ideal 
$\mf{m}\in {\M}(R)$, then $\alpha(X)\in {\E}_n(R[X])$. } \vp

Shortly after his proof of Serre's problem, Suslin-Kopeiko in \cite{SUSK}  
established an analogue of the Local-Global Pinciple for the 
elementary subgroup of the orthogonal group. Around the same time  
V.I. Kopeiko proved the analogous result for the  
elementary subgroup of the symplectic group. In this note we establish
an analogous Local-Global Principle for the elementary transvection subgroup 
of the automorphism group of projective, symplectic and orthogonal modules of 
global rank at least 1 and local rank at least 3. All previous work on this 
topic assumed that the global rank is at least 3. By definition the global 
rank or simply rank of a finitely generated projective $R$-module 
(resp. symplectic or orthogonal $R$-module) 
is the largest integer 
$r$ such that $\overset{r}\oplus R$ (resp. $\overset{r}\perp \mbb{H}(R)$) 
is a direct summand (resp. orthogonal summand) of the module. 
$\mbb{H}(R)$ denotes the hyperbolic plane. 

Using this Principle one can generalize well known facts regarding the group 
${\GL}_n(R)$ (${\Sp}_{2n}(R)$ or ${\O}_{2n}(R)$) of automorphisms 
of the free module $\overset{n}\oplus R$ of rank $n$ (free hyperbolic module 
$\overset{n}\perp \mbb{H}(R)$ of rank $n$)  to the automorphism group 
of finitely generated projective (symplectic or orthogonal) modules of 
global rank at least 1 and satisfying the local condition mentioned above. 
Specifically, we shall show that the elementary transvection subgroup 
is normal and the full automorphism group modulo its elementary transvection 
subgroup is nilpotent-by-abelian whenever stable dimension is finite. These 
generalize results by
Suslin and Kopeiko in \cite{KOP}, \cite{SUS}, \cite{SUSK}, Taddei in \cite{T}, 
first author in \cite{Bak}, Vavilov and Hazrat in \cite{HV}, and others.  
We treat the above three groups uniformly.

Our main results are as follows: \vp 

Let $Q$ denote a projective, symplectic or orthogonal module of global rank 
$\ge $ 1 and satisfying the local conditions stated above. Let
\begin{align*} 
{\G}(Q) & =  \tn{the automorphism group of $Q$}, \\
{\T}(Q) & =  \tn{the subgroup generated by transvections, and} \\
{\ET}(Q) & =  \tn{the subgroup generated by elementary transvections}.
\end{align*} 

{\bf Theorem 1.} {\it Let $R$ be a commutative ring with identity and let  
$\alpha(X)\in {\G}(Q[X])$, with $\alpha(0)={\I}_n$. If 
$\alpha_{\mf{m}}(X)\in {\ET}(Q_{\mf{m}}[X])$, for every maximal ideal 
$\mf{m}\in {\M}(R)$, then $\alpha(X)\in {\ET}(Q[X])$. } \vp 

{\bf Theorem 2.} {\it ${\T}(Q)={\ET}(Q)$. 
Hence ${\ET}(Q)$ is a normal subgroup of ${\G}(Q)$. } \vp 

By applying  Local-Global Principle (Theorem 1) we prove \vp

{\bf Theorem 3.} {\it
The factor group $\frac{{\G}(Q)}{{\ET}(Q)}$ is nilpotent-by-abealian, 
when the stable dimension $($i.e. Bass-Serre dimension$)$ is finite.} \vp 

To prove the result
we use the ideas of first author in \cite{Bak}, 
where he has shown that 
the group ${\GL}_n(R)/{\E}_n(R)$ is nilpotent-by-abelian for $n\ge 3$, 
but we avoid the functorial construction of the descending central series. 

\section{Preliminaries}

\begin{de}\tn{Let $R$ be an associative ring with identity.
The following condition was introduced by H. Bass:}

\tn{{\bf $(R_m)$} for every $(a_1,\dots, a_{m+1})
\in {\Um}_{m+1}(R)$, there are $\{x_i\}_{(1\le i\le m)}\in R$ such that} 
$(a_1+a_{m+1}x_1)R+\cdots+(a_m+a_{m+1}x_m)R=R$.

\tn{The condition $(R_m)\Ra (R_{m+1})$ for every $m> 0$. Moreover, for any 
$n\ge m+1$ the condition $(R_m)$ implies $(R_n)$ with $x_i=0$ for $i\ge m+1$.}

\tn{By {\bf stable range} for an associative ring $R$  
we mean the least $n$ such that $(R_n)$ holds.} 
\tn{The integer $n-1$ is called the {\bf stable 
dimension} of $R$ and is denoted by ${\sdim}(R)$.}
\end{de} 
\begin{lm} \label{sd1} $({\it cf.} \cite{B2})$
If $R$ is a commutative noetherian ring with identity of Krull dimension $d$, 
then ${\sdim}(R)\le d$.
\end{lm} 
\begin{de} \tn{A row vector $(a_1,\dots,a_n)\in R^n$ is said to be 
{\bf unimodular} in $R$ if $\sum_{i=1}^nRa_i=R$. The set of unimodular vectors 
of length $n$  in $R$ is denoted by $\tn{Um}_n(R)$. For an ideal $I$, 
$\tn{Um}_n(R,I)$ will denote the set of those unimodular vectors which are 
$(1,0,\dots,0)$ modulo $I$. }
\end{de}
\begin{de}\tn{Let $M$ be a finitely generated left module over a ring 
$R$. An element $m$ in $M$ is said to be {\bf unimodular} in $M$ if 
$Rm\cong R$  and $Rm$ is a direct summand of $M$, {\it i.e.} if there exists 
a finitely generated $R$-submodule $M'$ such that $M\cong Rm\oplus M'$.}
\end{de} 
\begin{de} \tn{For an element $m\in M$, one can attach an ideal, called the 
{\bf order ideal} of $m$ in $M$, {\it viz.} 
${\O}_M(m)=\{f(m)|f\in M^{*}=\tn{Hom}(M,R)\}$.
Clearly, $m$ is unimodular if and only if $Rm = R$ and ${\O}_M(m)=R$. }
\end{de}

\begin{de} \tn{Following H. Bass (\cite{B2}, pg. 167) we define 
a transvection of a finitely generated left $R$-module as follows: Let 
$M$ be a finitely generated left $R$-module. 
Let $q\in M$ and $\varphi\in M^{*}$ with $\varphi(q)=0$ . 
An automorphism of $M$ of the form $1+\varphi_q$ (defined by 
$\varphi_q(p)=\varphi(p)q$, for $p\in M$), will be called a 
{\bf transvection} of $M$ if either $q\in {\Um}(M)$ or $\varphi\in 
{\Um}(M^{*})$. We denote by Trans($M$)
the subgroup of Aut$(M)$ generated by transvections of $M$. }
\end{de}
\begin{de} \tn{Let $M$ be a finitely generated left $R$-module. 
The automorphisms of the form $(p,a)\mapsto (p+ax,a)$ and 
$(p,a)\mapsto (p,a+\psi(p))$, where $x\in M$ and $\psi
\in M^{*}$, are called {\bf elementary transvections} of $M\oplus R$. (It 
is easily verified that these automorphisms are transvections.)  
The subgroup  of ${\tran}(M\oplus R)$ generated by the
elementary transvections is denoted by ${\fl}(M\oplus R)$.}
\end{de} 
\begin{de} \tn{Let $R$ be an associative ring with identity.
To define other classical modules, we need an involutive antihomomorphism
({\bf involution}, in short) $*:R\ra R$ ({\it i.e.},  
$(x-y)^{*}=x^{*}-y^{*}$, $(xy)^{*}=y^{*}x^{*}$ and $(x^{*})^{*}=x$, for any 
$x,y\in R$). We assume that $1^{*}=1$. 
For any left $R$-module $M$ the involution induces a left module 
structure to the right $R$-module $M^{*}$ = Hom$(M,R)$ given by 
$(xf)v=(fv)x^{*}$, where $v\in M$, $x\in R$ and $f\in M^{*}$. 
Any right $R$ module can be viewed as a left $R$-module via the 
convention $ma=a^{*}m$ for $m\in M$ and $a\in R$. Hence if $M$ is a left 
$R$-module, then ${\O}_M(m)$ has a left $R$-module 
structure  with scalar multiplication given by $\lambda f(m)=f(\lambda m)$.}
\end{de} 
{\bf Blanket Assumption:} Let $A$ be an $R$-algebra, 
where $R$ is a commutative ring with identity, such that $A$ is finite as a 
left $R$-module.
Let $A$ possesses an involution 
$*:r\mapsto \bar{r}$, for $r\in A$. For a matrix $M=(m_{ij})$ over $A$ we  
define $\ol{M}=(\ol{m}_{ij})^t$. Let 
$\psi_1=\begin{pmatrix} 0 & 1 \\ -1 &0 \end{pmatrix}$, 
$\psi_n=\psi_{n-1}\perp \psi_1$ for $n>1$; and $\widetilde{\psi}_1=
\begin{pmatrix} 0 & 1 \\ 1 &0 \end{pmatrix}$, $\widetilde{\psi}_n=
\widetilde{\psi}_{n-1}\perp \widetilde{\psi}_1$, for $n>1$.
For a column vector $v\in A^n$ we write 
$\widetilde{v}=\bar{v}^t \psi_n$ in the symplectic case and 
$\widetilde{v}= \bar{v}^t \widetilde{\psi_n}$ in the orthogonal case. 
We define a form $\langle\, ,\rangle$ as follows: 
$$\langle v,w\rangle =\begin{cases} v^t\cdot w  & \mbox{in the linear case} \\
\widetilde{v}\cdot w & \mbox{otherwise.} \end{cases}$$
(Viewing $M$ as a right $A$-module we can assume the linearity).

Since $R$ is commutative, we can assume that 
the involution ``$*$'' defined on $A$ is trivial over 
$R$. We shall always assume that 2 is invertible in the ring $R$ while 
dealing with the symplectic and the orthogonal cases. 
\begin{de}\tn{A {\bf symplectic} ({\bf orthogonal}) $A$-module is a pair 
$(P,\langle \, ,\rangle)$, where 
$P$ is a projective left $A$-module of even rank and 
$\langle \, ,\rangle:P\times P\ra A$ is a 
non-singular ({\it i.e.} $P\cong P^{*}$ by 
$x\mapsto \langle x,\cdot \rangle$)  alternating (symmetric) bilinear 
form). }
\end{de}  
\begin{de} \tn{
Let $(P_1,{\langle \, ,\rangle}_1)$ and $(P_2,{\langle \, ,\rangle}_2)$ be two
symplectic (orthogonal) left $A$-modules. Their {\bf orthogonal sum} is the 
pair $(P,\langle \, ,\rangle)$, where $P=P_1\oplus P_2$ and the inner product 
is defined by $\langle(p_1,p_2),(q_1,q_2)\rangle={\langle p_1,q_1\rangle}_1 + 
{\langle p_2,q_2\rangle}_2$. Since this form is also non-singular we shall 
henceforth denote $(P,{\langle \, ,\rangle})$ by $P_1\perp P_2$ and called 
the orthogonal sum of $(P_1,{\langle \, ,\rangle}_1)$ and 
$(P_2,{\langle \, ,\rangle}_2)$ (if ${\langle \, ,\rangle}_1$ and 
${\langle \, ,\rangle}_2$ are clear from the context).}
\end{de} 
\begin{de}\tn{For a projective left $A$-module $P$ of rank $n$, we define 
$\mbb{H}(P)$ of global rank rank
$n$ supported by $P\oplus P^{*}$, with form $\langle (p,f),(p',f')\rangle
=f(p')-f'(p)$ for the symplectic modules and $f(p')+f'(p)$ for the 
orthogonal modules. 
There is a unique non-singular alternating (symmetric) bilinear form 
$\langle \, ,\rangle$ on the $A$-module $\mbb{H}(A) = A\oplus A^{*}$ 
(up to scalar multiplication 
by $A^{*}$) namely $\langle(a_1,b_1),(a_2,b_2)\rangle=a_1b_2-a_2b_1$ in 
the symplectic case and $a_1b_2+a_2b_1$ in the orthogonal case.}
\end{de} 
\begin{re}\tn{
A bilinear form $\langle \, ,\rangle$ induces a homomorphism $\Psi:P\ra P^{*} 
=\tn{Hom}(P,A)$, defined by $\Psi(p)(q)=\langle p,q\rangle$. 
The converse is also true since $2$ is invertible in $A$. If 
$\langle \, ,\rangle$ 
is symmetric, then one has $\Psi = \Psi^{*}$, and if $\langle \, ,\rangle$ is 
alternating, then one has $\Psi + \Psi^{*}=0$, under the canonical isomorphism 
$P\cong P^{**}$. }
\end{re} 
\begin{de} \tn{
An {\bf isometry} of a symplectic (orthogonal) module $(P,\langle \, 
,\rangle)$ is an automorphism of $P$ which fixes the bilinear form. 
The group of isometries of $(P, \langle \, ,\rangle)$ is denoted by 
${\Sp}(P)$ for the symplectic modules and 
${\O}(P)$ for the orthogonal modules. }
\end{de} 
\begin{de}\tn{Following Bass (\cite{B}) we define 
a symplectic transvection as follows: 
Let $\Psi:P\ra P^{*}$ be an induced isomorphism. 
Let $\alpha:A\ra P$ be a $A$-linear map defined by $\alpha(1)=u$. Then 
$\alpha^{*}\Psi\in P^{*}$ is defined by $\alpha^{*}\Psi(p)=
\langle u,p\rangle$. Let 
$v\in P$ be such that $\alpha^{*}\Psi(v)=\langle u,v\rangle=0$. An
automorphism $\sigma_{(u,v)}$ of $(P,\langle \, ,\rangle)$ of the form }
$$\sigma_{(u,v)}(p)=p+\langle u,p\rangle v+\langle v,p\rangle u+\langle u,
p\rangle u$$ 
\tn{for $u,v\in P$ with $\langle u,v\rangle=0$ will be called a
{\bf symplectic transvection} of $(P,\langle \, ,\rangle)$ if either 
$v\in {\Um}(P)$ or $\alpha^{*}\Psi\in {\Um}(P^{*})$. 
(Viewing $P$ as a right $A$-module we can assume the linearity).}
\end{de}

Since $\langle\sigma_{(u,v)}(p_1),\sigma_{(u,v)}(p_2)\rangle =\langle p_1,p_2
\rangle$, $\sigma_{(u,v)}\in {\Sp}(P,\langle \, ,\rangle)$.  
Note that $\sigma_{(u,v)}^{-1}(p)=p-\langle u,p\rangle v-\langle v,p\rangle u-
\langle u,p\rangle u$. 

The subgroup of ${\Sp}(P,\langle \, ,\rangle)$ generated by the symplectic 
transvections is denoted by ${\tran}_{\Sp}(P)$.
\begin{de}\tn{
The (symplectic) transvections of $(P\perp A^2)$ of the form 
$(p,b,a)\mapsto (p+aq,b-\langle p,q\rangle+a,a)$ and 
$(p,b,a)\mapsto (p+bq,b,a +\langle p,q\rangle -b)$, where
$a,b\in A$ and $p,q\in P$, are called 
{\bf elementary symplectic transvections}. 
The subgroup of ${\tran}_{\Sp}(P\perp A^2)$ 
generated by the elementary symplectic transvections is denoted by 
${\fl}_{\Sp}(P\perp A^2)$. }
\end{de}

In a similar manner we can define a transvection $\tau_{(u,v)}$ for an 
orthogonal module $(P, \langle ,\rangle)$. For this we need to assume 
that $u,v\in P$ are isotropic, {\it i.e.} 
$\langle u,u\rangle=\langle v,v\rangle=0$.  
\begin{de} \tn{An automorphism $\tau_{(u,v)}$ of $(P,\langle \, ,\rangle)$ 
of the form }
\end{de} 
$$\tau_{(u,v)}(p)=p-\langle u,p\rangle v+\langle v,p\rangle u$$ for 
$u,v\in P$ with $\langle u,v\rangle=0$ will be called an 
{\bf isotropic (orthogonal) transvection} of $(P,\langle \, ,\rangle)$ if 
either $v\in {\Um}(P)$ or $\alpha^{*}\Psi\in {\Um}(P^{*})$, 
(see \cite{HM}, pg. 214).

One checks that $\tau_{(u,v)}\in {\O}(P,\langle \, ,\rangle)$ and  
$\tau_{(u,v)}^{-1}(p)=p+\langle u,p\rangle v-\langle v,p\rangle u$.

The subgroup of ${\O}(P,\langle \, ,\rangle)$ generated by the isotropic 
orthogonal transvections is denoted by ${\tran}_{\O}(P)$. 
\begin{de} \tn{The isotropic orthogonal transvections of 
$(P\perp A^2)$ of the form 
$(p,b,a)\mapsto (p-aq,b+\langle p,q\rangle,a)$ and 
$(p,b,a)\mapsto (p-bq,b,a-\langle p,q\rangle)$, where
$a,b\in A$ and $p,q\in P$, are called 
{\bf elementary orthogonal transvections}. 
The subgroup of ${\tran}_{\O}(P\perp A^2)$ 
generated by elementary orthogonal transvections is denoted by 
${\fl}_{\O}(P\perp A^2)$. }
\end{de}
\begin{nt} \label{note}
\tn{In the sequel $P$ will denote either a finitely generated 
projective left $A$-module of rank 
$n$, a symplectic left $A$-module or an orthogonal left $A$-module of even rank 
$n=2r$ with a fixed form $\langle \, ,\rangle$. And 
$Q$ will denote $P\oplus A$ in the linear case, and $P\perp A^2$,
otherwise. To denote $(P\oplus A)[X]$ in the linear case and 
$(P\perp A^2)[X]$, otherwise, we will use the notation $Q[X]$. 
We assume that the rank of projective module $n\ge 2$, 
when dealing with the linear case, and $n\ge 6$, when considering the 
symplectic and the orthogonal cases. For a finitely generated projective 
$A$-module $M$ we use the notation }
\tn{${\G}(M)$ to denote ${\aut}(M)$, ${\Sp}(M,\langle \, ,\rangle )$ and 
${\O}(M,\langle \, ,\rangle )$ respectively; 
${\es}(M)$ to denote ${\SL}(M)$, ${\Sp}(M,\langle \, ,\rangle )$ and 
${\SO}(M,\langle \, ,\rangle )$ respectively; 
${\T}(M)$ to denote ${\tran}(M)$, ${\tran}_{\Sp}(M)$ and 
${\tran}_{\O}(M)$ respectively; and 
$\tn{ET}(M)$ to denote ${\fl}(M)$, ${\fl}_{\Sp}(M)$ and 
${\fl}_{\O}(M)$ respectively.}

\tn{
The reader should be able to easily verify that if $R$ is a reduced ring and 
$P$ a free $R$-module, i.e.
if $P = R^r$ (in the symplectic and the orthogonal cases we assume that $P$
is free with the standard bilinear form), then ETrans$(P) \supset$ E$_r(R)$,
ETransSp$(P) \supset $ ESp$(R)$ and ETrans$_O(P) \supset $ EO$_r(R)$, for $r
\geq 3$, in the linear case, and for $r \geq
6$, in the symplectic (and orthogonal) case.}

{\it Equality in all these cases will follow from Lemma \ref{free1} below.}  
\end{nt}  
We shall assume \vp\\
{\bf (H1) for every maximal ideal $\mf{m}$ of $A$, the symplectic (orthogonal) 
module $Q_{\mf{m}}$ is isomorphic 
to $A^{2n+2}_{\mf{m}}$ with the standard bilinear form $\mbb{H}(A_{\mf{m}}^{n+1})$.} 
\vp\\
{\bf (H2) for every non-nilpotent $s\in A$, if the projective module 
$Q_s$ is free $A_s$-module, then the symplectic (orthogonal) module 
$Q_s$ is isomorphic to $A^{2n+2}_s$ with the standard bilinear form 
$\mbb{H}(A_s^{n+1})$.}

\begin{re} \tn{Note that ${\T}(P)$ is a normal subgroup of ${\G}(P)$.
Indeed, for $\alpha\in {\G}(P)$ in the linear case we have 
$\alpha(1+\varphi_q)\alpha^{-1}
=1+(\varphi\alpha^{-1})_{\alpha(q)}$. In the symplectic case we can write 
$\sigma_{(u,v)}=1+\sigma'_{(u,v)}$ and similarly, 
$\alpha(1+\sigma'_{(u,v)})\alpha^{-1}
=1+(\sigma'\alpha^{-1})_{(\sigma'(u),\sigma'(v))}$. 
Similar argument will also holds for the orthogonal case.}
\end{re}

\begin{lm} \label{free1}
If the projective module $P$ is free of finite 
rank $n$ is free $($in the symplectic and the orthogonal 
cases we assume that the projective module is free with the standard bilinear 
form$)$, then 
${\tran}(P)={\E}_{n}(R)$, 
${\tran}_{\Sp}(P)= {\ESp}_{n}(R)$ and 
${\tran}_{\O}(P)= {\EO}_{n}(R)$ for $n\ge 3$ in the 
linear case and for $n\ge 6$ otherwise. 
\end{lm} 
{\bf Proof.}  In the linear case, for $p\in P$ and $\varphi\in P^{*}$ if 
$P=R^n$ then $\varphi_{p} : R^{n}\ra R\ra R^{n}$. Hence 
$1+\varphi_{p}=I_{n}+v.w^t$ for some row vector $v$ and column vector 
$w$ in $R^{n}$. Since $\varphi(p)=0$, it follows that 
$\langle v,w\rangle=0$. Since either $v$ or $w$ is unimodular, 
it follows that
$1+\varphi_{p}=I_{n}+v.w^t \in{\E}_{n}(R)$. 
Similarly, in the non-linear cases we have 
$\sigma_{(u,v)}(p)=I_{n}+v.\widetilde{w}+w.\widetilde{v}$, and  
$\tau_{(u,v)}(p)=I_{n}+v.\widetilde{w}-w.\widetilde{v}$, where 
either $v$ or $w$ is unimodular and $\langle v,w\rangle=0$. (Here 
$\sigma_{(u,v)}$ and $\tau_{(u,v)}$ are as in the definition of symplectic 
and orthogonal transvections.) Historically, these are known to be elementary 
matrices - for details see \cite{SUS} for the linear case, \cite{KOP} for the 
symplectic case, and \cite{SUSK} for the orthogonal case. \hb
\begin{re} \tn{Lemma \ref{free1} holds for $n=4$ in 
the symplectic and the orthogonal cases. This will follow from Remark 
\ref{free2}.}
\end{re}
\begin{re} \label{free2} 
\tn{${\ESp}_4(A)$ is a normal subgroup of ${\Sp}_4(A)$ by 
(\cite{KOP}, Corollary 1.11). Also 
${\ESp}_4(A[X])$ satisfies the Dilation Principle and the Local-Global 
Principle by (\cite{KOP}, Theorem 3.6). 
Since we were intent on a uniform proof, these cases 
have not been covered by us.}
\end{re} 
\begin{nt} \label{note1}
\tn{When $P=A^n$ ($n$ is even in the non-linear cases), we also use the notation
${\G}(n,A)$, ${\es}(n,A)$ and 
${\E}(n,A)$ for ${\G}(P)$, ${\es}(P)$ and ${\T}(P)$ respectively. 
We denote the usual standard elementary generators of ${\E}(n,A)$ by 
$ge_{ij}(x)$, $x\in A$. $e_i$ will denote the column vector  
$(0,\dots,1,\dots,0)^t$ (1 at the i-th position).}
\end{nt}
\begin{re} \label{endo}
\tn{Note that if $\alpha\in \tn{End}(Q)$ then 
$\alpha$ can be considered as a matrix of the form
{\small $\left( \begin{array}{cc} \tn{End}(P) & \tn{Hom}(P,A) \\ 
\tn{Hom}(A,P) & \tn{End}(A) \end{array}\right)$} in the linear case. 
In the non-linear cases one has a similar matrix for $\alpha$ of the form 
{\small $\left( \begin{array}{cc} \tn{End}(P) & \tn{Hom}(P,A\oplus A) \\ 
\tn{Hom}(A\oplus A,P) & \tn{End}(A\oplus A) \end{array}\right)$}}.
\end{re}
\section{\large Local-Global Principle for the Transvection Groups} 
~~~In this section we deduce an analogue of Quillen's Local-Global Principle 
for the linear, symplectic and isotropic orthogonal transvection groups. 
\begin{pr} \label{dil} {\bf (Dilation Principle)} 
Let $A$ be an associative $R$-algebra such that $A$ is finite as a left 
$R$-module and $R$ be a commutative ring with identity. 
Let $P$ and $Q$ be as in \ref{note}. Assume that $($H2$)$ holds.
Let $s$ be a non-nilpotent in $R$ such that $P_s$ is free, and 
let $\sigma(X)\in {\G}(Q[X])$ with $\sigma(0)=\tn{Id}$. Suppose 
$$\sigma_s(X)\in \begin{cases} {\E}(n+1,A_s[X]) & \mbox{in the linear case, }\\
{\E}(2n+2,A_s[X])& \mbox{ otherwise. } 
\end{cases}$$ 
Then there exists $\hat{\sigma}(X)\in {\rm ET}(Q[X])$ 
and $l>0$ such that $\hat{\sigma}(X)$ localizes to 
$\sigma(bX)$ for some $b\in (s^l)$ and $\hat{\sigma}(0)=\tn{Id}$.
\end{pr} 

First we state the following useful lemmas.
\begin{lm} \label{hom} 
Let $R$ be a ring and $M$ be a finitely presented left $($right$)$ 
$R$-module and $N$ be 
any $R$-module. Then we have a natural isomorphism: 
$$\gamma: \tn{Hom}_R(M,N)[X]\ra \tn{Hom}_{R[X]}(M[X],N[X]).$$ 
\end{lm} 
\begin{lm} \label{frac}
Let $S$ be a multiplicative closed subset of a ring $R$. Let $M$ 
be a finitely presented left $($right$)$ $R$-module and $N$ be any $R$-module. 
Then we have a natural isomorphism: 
$$\eta: S^{-1}(\tn{Hom}_R(M,N))\ra \tn{Hom}_{S^{-1}R}(S^{-1}M,S^{-1}N).$$
\end{lm} 
\begin{lm} 
\label{key4} \tn{(See \cite{BRK})} 
If $\eps=\eps_1\eps_2\cdots\eps_r$, 
where each $\eps_j$ is standard elementary generator, then for any $r>0$, 
and for any $(p,q)\in \mbb{N}\times\mbb{N}$,
$$\eps ge_{pq}(X^{2^rm}Y)\eps^{-1}=\underset{t=1}{\overset{k}\Pi}
ge_{p_tq_t}(X^m h_t(X,Y)),$$ for some $h_t(X,Y)\in R[X,Y]$, 
$(p_t,q_t)\in \mbb{N}\times\mbb{N}$, and for some $k>0$.
\end{lm} 
{\bf Proof of Proposition \ref{dil}.} Since elementary transvections can 
always be lifted, we can and hence assume that $R$ is reduced. 
We show that there exists 
$l>0$ such that $\sigma(bX)\in {\rm ET}(Q[X])$ for all $b\in (s^k)R$, for all $k\ge l$. 

As $\sigma(0)=$ Id, we can write $\sigma_s(X)=\underset{k}\Pi
\gamma_k ge_{i_k j_k}(X\lambda_k(X))\gamma_k^{-1}$, 
where $\lambda_k(X)\in A_s[X]$. Hence for $d> 0$,
$\sigma_s(XT^{2d})=\underset{k}\Pi\gamma_k ge_{i_k j_k}(XT^{2d}
\lambda_k(XT^{2d}))\gamma_k^{-1}$, for some $\gamma_k$ in ${\E}(n+1,R_s)$ 
in the linear case, and in ${\E}(2n+n,R_s)$ in the non-linear cases. Let 
$\gamma_k=\eps_1\eps_2\ldots\eps_r$, and $d=2^{r-1}$. Now apply Lemma 
\ref{key4} for $m=1$, $X=T$, and $Y=1$. Then use the fact that $i\ne 1\ne j$, 
$$ge_{ij}(T^2\mu(X))=[ge_{i1}(T\mu(X)), ge_{1j}(T)],$$ in the linear case, and 
$$ge_{ij}(T^2\mu (X))=\begin{cases}
[ge_{i1}(T\mu (X)), ge_{1j}(T)]  & \mbox{ if }  i\ne \sigma(j)  \\
ge_{1j}(T\mu (X))[ge_{1l}(-T\mu (X)), ge_{i1}(-T)] & \mbox{ if } i=\sigma (j)\\
[ge_{\sigma(i)1}(T\mu (X)), se_{1\sigma(j)}(T)] & \mbox{ if } i\ne \sigma (j)\\
ge_{1\sigma{j}}(T\mu (X))[ge_{1\sigma(l)}(-T\mu (X)), ge_{\sigma(i)1}(-T)] & \mbox{ if } 
i=\sigma (j)
\end{cases}$$
in the non-linear cases for some $l\le n$, when $i+1$ is even and when 
$\sigma(i)+1$ is even respectively. 
Then for $d\gg 0$ we get 
$\sigma_s(XT^{2d})=\underset{t}\Pi ge_{p_t q_t}(T\mu_t(X))$,
for some $\mu_t(X)\in A_s[X]$ with $p_t=1$ or $q_t=1$.

Since $P_s$ is free $A_s$-module, 
\begin{align*} 
P_s[X,T] & \cong A_s^{n}[X,T]\cong P_s[X,T]^{*}    \tn{ ~in the linear case},\\
P_s[X,T] & \cong A_s^{2n}[X,T]\cong P_s[X,T]^{*}  \tn{ in the non-linear cases}.
\end{align*} 
Thus using the isomorphism, 
polynomials in $P_s[X,T]$ can be regarded as linear forms
which acts as follows: For $x=(x_1,\dots, x_k), y=(y_1,\dots, y_k)
\in A_s^k[X,T]$ ($k=n$ in the linear case and $k=2n$ in the symplectic case), 
$$\langle x,y\rangle = 
\begin{cases} xy^t & \mbox{ in the linear case, }\\
x\psi_{n} y^t & \mbox{ in the symplectic case, }\\
x\widetilde{\psi_{n}} y^t & \mbox{ in the orthogonal case. }\\
\end{cases}$$
\noindent (where $\psi_n$ denote the alternating matrix corresponding 
to the standard symplectic form
$\underset{i=1}{\overset{2n}\sum} e_{2i-1,i}-
\underset{i=1}{\overset{2n}\sum} e_{2i,2i-1}$ 
and $\tilde{\psi}_n$ denote the symmetric matrix corresponding to the 
standard hyperbolic form 
$\underset{i=1}{\overset{2n}\sum} e_{2i-1,i}+ 
\underset{i=1}{\overset{2n}\sum} e_{2i,2i-1}$).

First we consider the case when $p_t=1$. 
Let $p_1^{*},\dots,p_k^{*}$ be the standard basis of $P_s$ and let 
$s^mp_i^{*}\in P$ for some $m>0$ and $i=1,\dots, k$. Let $e_i^{*}$ be the 
standard basis of $A_s$. For $q_t=i$, consider the element 
$T\mu_t(X)e_i^{*}\in A_s^k[X,T]$ as an element in $P_s[X,T]^{*}$. Using 
Lemma \ref{hom}, we may say $T\mu_t(X)e_i^{*}$ is actually a polynomial 
in $T$. By Lemma \ref{frac}, there exists $k_1>0$ such that 
$k_1$ is the maximum power of $s$ occurring in the denominator of 
$\mu_t(X)e_i^{*}$. Choose $l_1\ge \tn{max}(k_1,m)$. 

Next suppose $q_t=1$. For $p_t=j$, $T\mu_t(X)e_j^{*}\in P_s[X,T]$. From 
Lemma \ref{frac} it follows that we can choose $k_2 >0$ such that $k_2$ is
the maximum power of $s$ occurring in $\mu_t(X)e_j^{*}$. Again, using Lemma 
\ref{hom}, we can regard $T\mu_t(X)e_j^{*}$ as a polynomial in $T$. Choose
$l_2\ge \tn{max}(k_2,m)$ and $l\ge \tn{max}(l_1,l_2)$. Now applying 
homomorphism $T\mapsto s^lT$ it follows that $\sigma(bXT^{2d})$ is defined 
over $Q[X]$. Putting $T=1$, by usual Dilation Principle there exists 
$l>0$ such that $\hat{\sigma}(X)\in {\rm ET}(Q[X])$ localizes to 
$\sigma(bX)$ for some $b\in (s^l)$ and $\hat{\sigma}(0)=\tn{ Id}$. \hb
\begin{lm} \label{sol2} 
Let $A$ be an associative $R$-algebra such that $A$ is finite as a left 
$R$-module and $R$ be a commutative ring with identity.
Let $\alpha\in {\es}(n,R)$ and $I$ be an ideal 
contained in the nil radical $\tn{Nil}(R)$ of $R$. 
Let `bar' denote the reduction 
modulo $I$. If $\ol{\alpha}\in {\E}(n,\ol{A})$, then $\alpha\in {\E}(n,A)$.
\end{lm} 
{\bf Proof.} This is easy to verify. (The proof of Lemma 
\ref{sol3a} later is similar; and its argument can be used to 
prove this lemma). \hb
\begin{tr} {\bf (Local-Global Principle)} \label{lgt}
Let $A$ be an associative $R$-algebra such that $A$ is finite as a left 
$R$-module and $R$ be a commutative ring with identity. 
Let $P$ and $Q$ be as in \ref{note}. Assume that $($H1$)$ holds.
Suppose $\sigma(X)\in {\G}(Q[X])$ with
$\sigma(0)=\tn{Id}$. If 
$$\sigma_{\mf{p}}(X)\in \begin{cases} {\E}(n+1,A_{\mf{p}}[X]) & 
\mbox{in the linear case, }\\
{\E}(2n+2,A_{\mf{p}}[X])& \mbox{ otherwise } 
\end{cases}$$ 
for all $\mf{p}\in {\sp}(R)$, then $\sigma(X)\in {\rm ET}(Q[X])$.
\end{tr} 
{\bf Proof.} Follows from the similar argument as in the proof of 
$(4)\Ra (3)$ in Theorem 3.1 of \cite{BRK}. \hb
\begin{co} \label{new} Assume that $($H2$)$ holds.
Let $\tau(X)\in {\G}(Q[X])$, with 
$\tau(0)=\tn{Id}$. If $\tau_s(X)\in {\ET}(Q_s[X])$, 
and $\tau_t(X)\in {\ET}(Q_t[X])$, for some $s$, $t\in R$ with that $Rs+Rt=R$, 
then $\tau(X)\in {\ET}(Q[X])$. 
\end{co}
\begin{lm} \label{mat}
Let $R$ be a commutative ring with identity. If 
$\alpha=(a_{ij})$ is an $r\times r$ matrix over $R$ with all entries 
nilpotent, then $\alpha$ is nilpotent. 
\end{lm} 
{\bf Proof.} Let $a_{ij}^l = 0$, for all $i,j\in \{1,\dots,r\}$. 
Now $\alpha^2$ has entries which are 
homogeneous polynomial of degree $2$ in the $a_{ij}$'s. Consequently, 
$\alpha^4$ has entries which are homogeneous polynomial of degree $4$ in the 
$a_{ij}$'s, and so on.  Therefore, $\alpha^{2^m} = 0$, if $2^m > l r^2$, 
by the Pigeon Hole Principle. \hb
\begin{co} \label{radical}
Let $A$ be an associative $R$-algebra such that $A$ is finite as a left 
$R$-module and $R$ be a commutative ring with identity. Let 
`bar' denote the reduction modulo $\tn{Nil}(R)$. Assume that $($H1$)$ 
holds. Then for $\tau\in {\G}(Q)$, $\ol{\tau}\in \tn{ET}(\ol{Q})$ if and only 
if $\tau\in \tn{ET}(Q)$. 
\end{co} 
{\bf Proof.} First we show that for $\sigma (X)\in {\G}(Q[X])$, 
$\ol{\sigma(X)}\in \tn{ET}(\ol{Q[X]})\equi
\sigma(X)\in \tn{ET}(Q[X])$.  
Suppose $\ol{\sigma(X)}\in \tn{ET}(\ol{Q[X]})$. Then 
$\ol{\sigma_{\mf{p}}(X)}\in \tn{ET}(\ol{A_{\mf{p}}}[X])=
{\E}(n+1,\ol{A_{\mf{p}}}[X])$ and ${\E}(2n+2,\ol{A_{\mf{p}}}[X])$ in the linear 
and the non-linear cases respectively, for all 
$\mf{p}\in {\sp}(R)$. But then from Lemma \ref{sol2} it follows that 
$$\sigma_{\mf{p}}(X)\in \begin{cases} {\E}(n+1,{A_{\mf{p}}}[X]) 
\mbox{  in the linear case, and }\\
{\E}(2n+2,{A_{\mf{p}}}[X]) \mbox{  otherwise. }
\end{cases}$$ 
Hence by Theorem \ref{lgt}, $\sigma(X)\in \tn{ET} (Q[X])$. 

Now modifying $\tau$ by some $\eps \in \tn{ET}(Q)$ we assume that 
$\tau=\tn{Id}+\gamma$, where $\gamma\equiv 0$ modulo $\tn{Nil}(R)$. 
As the nilpotent entries are in $R$, which is a commutative ring, by 
Lemma \ref{mat}, $\gamma$ is nilpotent. 
Define $\theta(X)=\tn{Id}+X\gamma$. 
As $\ol{\theta(X)}=\tn{Id}\in \tn{ET}(\ol{Q[X]})$, from above it follows that
$\theta(X)\in \tn{ET}(Q[X])$. Whence $\tau=\theta(1)\in 
\tn{ET}(Q)$; as required. \hb 
\begin{tr} \label{tf} 
Let $A$ be an associative $R$-algebra such that $A$ is finite as a left 
$R$-module and $R$ be a commutative ring with identity. Let
$Q$ be as in \ref{note}. Assume that $($H1$)$ holds. 
Then ${\T}(Q)=\tn{ET}(Q)$. 
\end{tr}
{\bf Proof.} Using Corollary \ref{radical} we assume that $A$ is reduced.  
By definition $\tn{ET}(Q)\subset {\T}(Q)$. To prove the converse 
assume $\tau\in {\T}(Q)$. Then there exists 
$\sigma(X)\in {\T}(Q[X])$ such that $\sigma(0)={\rm Id}$ and $\sigma(1)=\tau$. 
Now, for every ${\mf{p}}\in {\sp}(R)$,  
$$\sigma_{\mf{p}} \in \begin{cases} 
{\E}(n+1,A_{\mf{p}}[X]),  \mbox{  in the linear case, and } \\ 
{\E}(2n+2,A_{\mf{p}}[X]), \mbox{ otherwise.}
\end{cases}$$
 Therefore, by Theorem 
\ref{lgt} it follows that
$\sigma(X)\in \tn{ET}(Q[X])$. Whence $\tau=\sigma(1)\in \tn{ET}(Q)$; 
as required. \hb \vp

The next lemma was proved in the linear case in (\cite{BMR}, Proposition 4.1). 
They do not assume the existence of a unimodular element in $Q$ though, 
and only get a unipotent lift. In (\cite{BR}, Lemma 2.1) it is mentioned that 
if $Q$ has a unimodular element then the lift is a transvection.

\begin{co} \label{roy} Let $A$ be an associative $R$-algebra such that $A$ is 
finite as a left $R$-module and $R$ be a commutative ring with identity.
Let $Q$ be as in \ref{note}. If $I$ in an ideal in $R$, then the map 
${\T}(Q)\ra {\T}(Q/IQ)$ is surjective.
\end{co} 
{\bf Proof.} By Proposition \ref{tf}, ${\T}(Q/IQ)=\tn{ET}(Q/IQ)$. Since 
an elementary transvection can always be lifted to an elementary transvection, 
the result follows. \hb 

As a consequence of Theorem \ref{dil} and Theorem \ref{tf}, 
following L.N. Vaserstein's proof of Serre's conjecture 
(see \cite{LAM}, Chapter III, \S 2) we deduce the following Local-Global 
Principle for the action of the elementary subgroups of an extended projective, symplectic and orthogonal modules.
module.    
\begin{tr} \label{lift} 
Let $A$ be an associative $R$-algebra such that $A$ is finite as a left 
$R$-module and $R$ be a commutative ring with identity. 
Let $Q$ be as in \ref{note}. Assume that $($H1$)$ and $(H2)$ holds.
Let $v(X)=(p(X),a(X))$ be a unimodular element in $Q[X]$ with $v(0)
=(0,1)$. Suppose for all prime ideal $\mf{p}\in {\sp}(R)$ there exists 
$$\tau_{\mf{p}}(X)\in  {\T}(Q_{\mf{p}}[X])= \begin{cases} 
{\E}(n+1,A_{\mf{p}}[X])  \mbox{ in the linear case, and} \\
{\E}(2n+2,A_{\mf{p}}[X]) \mbox{ in non-linear cases} \end{cases}$$
with $\tau_{\mf{p}}(0)=I_n$ such that $v(X)\tau(X)=v(0)$. 
Then there exists $\tau(X)\in \tn{ET}(Q[X])$ with $\tau(0)=\tn{Id}$ 
such that $v(X)\tau(X)=v(0)$.
\end{tr}  

\section{\large The unstable ${\k}$-groups 
$\frac{{\es}(Q)}{{\ET}(Q)}$ are nilpotent}

~~~~In this section earlier results on unstable ${\k}$-groups of 
classical groups of A. Bak, R. Hazrat, and H. Vavilov have been 
uniformly generalized to classical modules. 

We prove Theorem 3 mentioned in the introduction. 
Before that we give a brief historical sketch about our result. 
Throughout this section we assume $R$ is a commutative ring with identity.

In \cite{Bak}, A. Bak defines a functorial filtration 
${\GL}_n (R)={\es}^{-1}{\rm L}_n(R)\supset {\es}^0{\rm L}_n(R) \lb \supset 
\cdots \supset {\es}^i{\rm L}_n (R) \supset \cdots \supset {\E}_n (R)$ 
of the general linear group ${\GL}_n (R)$, where $R$ is an associative ring 
with identity and $n \geq 3$, which is a descending central
series. His construction has its own merits; which we do 
not study here though, other than the fact that the quotient 
${\GL}_n(R)/{\E}_n (R)$ is nilpotent-by-abelian. A. Bak uses a 
localization-completion method; we show that the localization part 
suffices to get the desired result.

In \cite{HV}, R. Hazrat and N. Vavilov have shown: Let $\Phi$ be a reduced
irreducible root system of rank $\geq 2$ and $R$ be a commutative ring
such that its Bass-Serre dimension $\delta(R)$ is finite. Then for any   
Chevalley group ${\G}(\Phi, R)$ of type $\Phi$ over $R$ the quotient 
${\G}(\Phi, R)/{\E}(\Phi, R)$ is nilpotent-by-abelian. In particular, 
${\k}(\Phi, R)$ is nilpotent of class at most $\delta(R) + 1$. They use the
localization-completion method of A. Bak in \cite{Bak}, 
who showed that ${\k}(n,R)$ is nilpotent-by-abelian. Their main result is to 
construct a descending central series in the Chevalley group,
indexed by the Bass-Serre dimension of the factor-rings of the ground ring.
Our approach show that for classical groups the localization part suffices.

The precise statement of our theorem is the following: 

\begin{tr} \label{nilpo}
Assume the notation in \ref{note}. We assume that 
$(H1)$  and $(H2)$ holds and $R$ is noetherian. 
Let $d=\dim \,(R)$ and $t=\tn{local rank of } Q$.

The quotient 
group $Q/\T(Q)$ is nilpotent of class at most \tn{max} 
$(1, d+3-t)$ in the linear case and \tn{max}$(1, d+3-t/2)$ otherwise. 
\end{tr} 

(We assume $Q$ that global rank of $Q$ is at least 1 and lonal rank of $Q$ 
is at least 3 in all the above cases).

\begin{co} Let $d=\dim \,(R)$ and $t=\tn{local rank of } Q$. The quotient 
group $Q/\ET(Q)$
is nilpotent of class at most \tn{max} 
$(1, d+3-t)$ in the linear case and \tn{max}$(1, d+3-t/2)$ otherwise. 
\end{co}

Recall
\begin{de} \tn{ Let $H$ be a group. Define $Z^0=H$, $Z^1=[H,H]$ and 
$Z^i=[H,Z^{i-1}]$. Then $H$ is said to be nilpotent if $Z^r=\{e\}$ for
some $r>0$, where $e$ denotes the identity element of $H$. }
\end{de}
\begin{re} \tn{ The standard 
generators of ${\E}(n,R)$ satisfy the following relation:}
$$[ge_{ik}(x),ge_{kj}(y)]=ge_{ij}(zxy)$$
\tn{for  $x,y\in R$ and some $z\in \mbb{N}$ (fixed for the group), and 
$1\le i\ne j\ne k\le n$, $j\ne\sigma i, \sigma j$, where $\sigma$ is the
  permutation given by $2l\mapsto 2l-1$ and $2l-1\mapsto 2l$.}
\end{re} 
\begin{nt} \tn{
Let ${\es}(n,s^lR)$ be the subgroup of ${\es}(n,R)$ consisting 
of matrices which are identity modulo $s^lR$ and ${\es}(Q,s^lR)$ the 
subgroup of ${\es}(Q)$ consisting of the automorphisms with determinant 
$1$ which are identity modulo $s^lR$.}
\end{nt}

Let $J(R)$ denote the Jacobson radical of $R$. 
\begin{lm} \label{sol3a}
Let $\beta\in {\es}(n,R)$, with $\beta\equiv I_n$ modulo $I$, 
where $I$ is contained in the Jacobson 
radical $J(R)$ of $R$. Then there 
exists $\eps\in {\E}(n,R)$ such that $\beta\eps$= the diagonal matrix 
$[d_1,d_2,\dots,d_n]$, where each $d_i$ is a unit in $R$ and each 
$d_i$ satisfies $d_{\sigma i}{d}_i^{*} = 1$, and $\eps$ a product of 
elementary generators with each congruent to identity modulo $I$.
\end{lm} 
{\bf Proof.} The diagonal elements are units. Using this one can establish 
the result easily in the linear case.
In the symplectic and the orthogonal cases by multiplying 
from the right side by suitable standard 
elementary generators each of which is 
congruent to identity modulo $I$ we can 
make all the $(1,j)$-th entries zero for $j=3,\dots,n$. Since char $R\ne 2$, 
in the orthogonal case the $(1,2)$-th entry will be then automatically zero. 
In the symplectic case again by multiplying from right side by  suitable 
elementary generators we can make the $(1,2)$-th entry zero. Similarly, 
multiplying by left we can make the first column of $\beta$ to be 
$(d_1,0,\dots,0)^t$, where $d_1$ is a unit in $R$ and $d_1\equiv 1$ modulo 
$I$. Repeating the above process we can reduce the size of $\beta$. 
Note that after modifying the first row and the first column in the 
symplectic and 
the orthogonal cases the second row and column will automatically become 
$(0,d_2,0,\dots,0)^t$ for some unit $d_2$ in $R$, with $d_2\equiv 1$ modulo 
$I$. Repeating the process we can modify $\beta$ to the required form. \hb
\vp \\
{\bf Blanket assumption:} 
Henceforth we shall assume that the matrices have size at least $3\times 3$ 
when dealing with the linear case and at least $6\times 6$ when dealing 
with the symplectic and the orthogonal cases. 
\begin{lm} \label{sol3}
Let $R$ be a commutative ring and $s$ be a non-zero divisor in $R$.
Let $D$ denote the diagonal matrix $[d_1,\dots,d_n]$, where $d_i$ should be 
units and satisfy $d_{\sigma i}{d}_i^{*} = 1$, and 
$d_i\equiv 1$ modulo $(s^l)$ for $l\ge 2$. Then 
$$\left[ge_{ij}\left(\frac{a}{s} X \right), D\right]\subset 
{\E}(n,R[X])\cap {\es}(n,(s^{l-1})R).$$
\end{lm} 
{\bf Proof.} Let $d=d_id_j^{-1}$. Then 
$\left[ge_{ij}\left(\frac{a}{s} X \right), D\right]= 
ge_{ij}\left(\frac{a}{s} X\right) ge_{ij}\left(-\frac{a}{s} dX\right)$. Since 
$d_i, d_j\equiv 1$ modulo $(s^l)$ for $l\ge 2$, we can write $d=1+s^m\lambda$ 
for some $m>2$ and $\lambda\in R$. Hence 
{\small \begin{align*}
ge_{ij}\left(\frac{a}{s} X \right)ge_{ij}\left(-\frac{a}{s} dX\right) & =
ge_{ij}\left(\frac{a}{s} X\right) ge_{ij}\left(-\frac{a}{s} X\right)
ge_{ij}\left(-\frac{a}{s}s^m\lambda X\right)\\
& =  ge_{ij}\left(-\frac{a}{s}s^m\lambda X\right)\in 
{\E}(n,R[X])\cap {\es}(n,(s^{m-1})R).
\end{align*}} \hb 
\begin{lm} \label{sol4}
Let $R$ be a ring, $s\in R$ a non-zero divisor in $R$ and $a\in R$. 
Then for $l\ge 2$  
$$\left[ge_{ij}\left(\frac{a}{s} X\right), {\es}(n,s^lR)\right]
\subset {\E}(n,R[X]).$$
More generally, 
$\left[\eps(X), {\es}(n,s^lR[X])\right]\subset {\E}(n,R[X])$ for $l\gg 0$ and 
$\eps(X)\in {\E}(n,R_s[X])$. 
\end{lm} 
{\bf Proof.} First fix $(i,j)$ for $i\ne j$. 
Let $\alpha(X)=[e_{ij}\left(\frac{a}{s} X\right), \beta]$ for 
some $\beta\in {\es}(n,s^lR)$. 
As $l\ge 2$, it follows that 
$\alpha(X)\in {\es}(n,R[X])$. Since ${\E}(n,R[X])$ is a normal subgroup of 
${\es}(n,R[X])$, we get $\alpha_s(X)\in {\E}(n,R_s[X])$. 
Let  $B=1+sR$. We show that $\alpha_B(X)\in {\E}(n,R_B[X])$. Since 
$s\in J(R_B)$,  it follows from Lemma \ref{sol3a}
that we can decompose $\beta_B=\eps_1\cdots\eps_tD,$ where 
$\eps_i=ge_{p_iq_i}(s^l \lambda_i) \in {\E}(n,R_B)$; $\lambda_i\in R_B$
and $D$ = the diagonal matrix $[d_1,\dots,d_n]$ with $d_i$ is a unit in $R$ 
and $d_i\equiv 1$ modulo $(s^l)$ for $l\ge 2$; $i=1,\dots,n$. 
If $t=1$, then using the commutator law and Lemma \ref{sol3}
it follows that $\alpha_B(X)\in {\E}(n,R_B[X])$. Suppose $t>1$. Then 
\begin{align*}
\alpha_B(X) &
 =\left[ge_{ij}\left(\frac{a}{s} X\right), \eps_1 \cdots \eps_t D\right] \\
&  =\left[ge_{ij}\left(\frac{a}{s} X\right), \eps_1\right] \eps_1 
\left[ge_{ij}\left(\frac{a}{s} X\right),\eps_2 \cdots \eps_t D\right] 
\eps_1^{-1}
\end{align*}
and by induction each term is in ${\E}(n,R_B[X])$, hence 
$\alpha_B(X)\!\in \!{\E}(n,R_B[X])$. Since $\alpha(0)=I_n$, 
by the Local-Global Principle for the classical groups it follows that 
$\alpha(X)\in {\E}(n,R[X])$. \hb
\begin{co} \label{sol5}
Let $R$ be a ring, $s\in R$ be a non-zero divisor  in $R$ and $a\in R$. 
Then for $l\ge 2$  
$$\left[ge_{ij}\left(\frac{a}{s} \right), {\es}(n,s^lR)\right]
\subset {\E}(n,R).$$
More generally, 
$\left[\eps, {\es}(n,s^lR)\right]\subset {\E}(n,R)$ for $l\gg 0$ and 
$\eps\in {\E}(n,R_s)$. 
\end{co}  
\begin{lm} \label{tr5} 
Fix the notation as in \ref{note}. 
Let $s$ be a non-zero divisor in $R$ such that $P_s$ is free. 
Assume that $(H2)$ holds. Suppose $\tau \in {\T}(Q_s)$. 
Then for $l\gg 0$, $[\tau, {\es}(Q,s^lR)]\subset {\T}(Q)$.
\end{lm}  
{\bf Proof.} Let $\eta\in {\es}(Q,s^lR)$ and 
$\widetilde{\tau}(X)\in {\T}(Q_s[X])$ 
be an isotopy between the identity map and $\tau$; {\it i.e.} 
$\widetilde{\tau}(0)=$Id and 
$\widetilde{\tau}(1)=\tau$. Let $\alpha(X)=[\widetilde{\tau}(X),\eta]$. 
Now, since $\eta\equiv $ Id modulo $(s^l)$, $\eta=\tn{Id}+s^l\psi$ 
for some $\psi \in  \tn{End}(Q)$. Therefore, 
$\psi$ can be considered as a matrix as in Remark \ref{endo}. 
Hence $\widetilde{\tau}(X) \, \eta \, \widetilde{\tau}(X)^{-1}=
\tn{Id}+s^l\widetilde{\tau}(X)\psi\widetilde{\tau}(X)^{-1} \in {\es}(Q[X])$
for $l\gg 0$. As `${\T}$' is a normal subgroup of `${\es}$', 
it follows that $\alpha_s(X)\in {\T}(Q_s[X])$.  Let $B=1+sR$.
We show that $\alpha_B(X)\in {\T}(Q_B[X])$. 
Note that $s\in \tn{Jac}(R_B)$. Hence for all $\mf{m}\in \tn{Max} (R_B)$, 
$$(\eta_B)_{\mf{m}}\in {\es}((Q_B)_{\mf{m}},
s^l(R_B)_{\mf{m}})
=\begin{cases} {\E}({n+1},(R_B)_{\mf{m}}) & \mbox{ in the linear case }\\
{\E}({2n+2},(R_B)_{\mf{m}}) & \mbox{ otherwise. }\end{cases}$$
Therefore, by Lemma \ref{sol3}, 
$(\eta_B)_{\mf{m}}$ can be expressed as a product of elementary matrices over 
$(R_B)_{\mf{m}}$ with each being identity modulo $(s^l)$, and a diagonal
matrix $D=[d_1,\ldots, d_t]$, where $t=r+1$ in the linear case and $r+2$
otherwise, and $d_i$ are units in $R$ for $i=1,\ldots,t$. Let 
$(\eta_B)_{\mf{m}}=\Pi_{i=1}^k \eps_i D$, where 
$\eps_i$ is in ${\E}_{n+1}((R_B)_{\mf{m}})$ in the linear case and  in 
${\E}_{2n+2}((R_B)_{\mf{m}})$ otherwise, and $\eps_i=$Id modulo $(s^l)$. 
So, $(\alpha_B)_{\mf{m}}(X)=[\widetilde{\tau}(X),\eps_1\cdots\eps_k D].$
Hence by Lemma \ref{sol3} and Lemma \ref{sol4}, we get 
$$(\alpha_B)_{\mf{m}}(X)\in \begin{cases} {\E}({n+1},(R_B)_{\mf{m}}[X]) & 
\mbox{ in the linear case } \\
{\E}({2n+2},(R_B)_{\mf{m}}[X]) & \mbox{ otherwise. } \end{cases}$$ 
Hence by the L-G Principle for the tranvection 
groups we get $\alpha_B(X)\in {\T}(Q_B[X])$.
Therefore, it follows from Corollary \ref{new} that $\alpha(X)\in 
{\T}(Q[X])$. In particular, $[\tau,\eta]\in {\T}(Q)$. $\tn{ }$ \hb \vp \\
{\bf Proof of Theorem \ref{nilpo}.} Using Corollary \ref{radical} we may and do 
assume that $R$ is a reduced 
ring. Note that if $t\ge d+3$, then the group ${\es}(Q)/{\T}(Q)={\k}(Q)$, 
which is abelian and hence nilpotent. So we consider the case $t\le d+3$. 
Let us first fix a $t$. We prove the theorem by induction 
on $d=\dim R$. Let $H={\es}(Q)/{\T}(Q)$. Let $m=d+3-t$ and 
$\alpha=[\beta, \gamma]$ for some $\beta\in H$ and $\gamma\in Z^{m-1}$. 
Clearly, the result is true for $d=0$. 
Let $\widetilde{\beta}$ be the pre-image of $\beta$ under the map 
${\es}(Q)\ra {\es}(Q)/{\T}(Q)$.
Choose a non-nilpotent element $s$ in $R$ such that $P_s$ is free and 
$\widetilde{\beta}_s\in {\E}(n,A_s)$.
We define $\ol{H}={\es}(\ol{Q})/{\T}(\ol{Q})$,
where bar denote reduction modulo $s^l$ for some $l\gg 0$. By the  
induction hypothesis $\ol{\gamma}=\{1\}$ in $\ol{{\G}}(Q)$. 
Since ${\T}(Q)$ is a normal subgroup
of ${\es}(Q)$ for $n\ge 3$ in the linear case and for $n\ge 4$ otherwise,
by modifying $\gamma$ we may assume that
$\widetilde{\gamma}\in {\es}(Q,s^lA)$, where $\widetilde{\gamma}$ is the pre 
image of $\gamma$ in ${\es}(Q)$. Now by Lemma \ref{tr5} it follows that 
$[\widetilde{\beta},\widetilde{\gamma}]\in {\T}(Q)$. 
Hence $\alpha=\{1\}$ in $H$. \hb \vp 

{\bf Acknowledgement:} The authors thank W. van der Kallen for pointing out 
that the argument here does not work with sdim instead of dim in Theorem 4.1, and Corollary 4.2.

\addcontentsline{toc}{chapter}{Bibliography} 
{\it Department of Mathematics, University of Bielefeld,
Germany. }\\
{\it email: bak@mathematik.uni-bielefeld.de }\vp\\
{\it Indian Institute of Science Education and Research, Kolkata, India.}\\
{\it email: rabeya.basu@gmail.com, rbasu@iiserkol.ac.in} \vp \\
{\it Tata Institute of Fundamental Research, 
Mumbai, India.}\\
{\it email: ravi@math.tifr.res.in}

\end{document}